\documentclass[11pt]{amsart}

\usepackage[a4paper,margin=28mm]{geometry}
\usepackage{amsmath,amssymb,amsfonts,amsthm,mathrsfs}
\usepackage{mathtools}
\usepackage{microtype}
\usepackage{hyperref}
\usepackage{xcolor}
\hypersetup{
    colorlinks,
    linkcolor={red!50!black},
    citecolor={blue!100!black},
    urlcolor={blue!100!black}
}

\numberwithin{equation}{section}
\newtheorem{theorem}{Theorem}[section]
\newtheorem{proposition}[theorem]{Proposition}
\newtheorem{lemma}[theorem]{Lemma}

\theoremstyle{remark}
\newtheorem*{remark}{Remark}
\newcommand{\T}{\mathbb T}
\newcommand{\R}{\mathbb R}
\newcommand{\Z}{\mathbb Z}
\newcommand{\calL}{\mathcal L}
\newcommand{\calH}{\mathcal H}
\newcommand{\Maj}{\mathfrak M}
\newcommand{\Min}{\mathfrak m}
\newcommand{\eps}{\varepsilon}

\title[Liouville exponential sums]
{On the lower bound for the $L^1$ norm of Liouville exponential sums. II.}

\author[B. Krause]{Ben Krause}
\address{School of Mathematics, University of Bristol, Fry Building,
Woodland Road, Bristol BS8 1UG, UK}
\email{ben.krause@bristol.ac.uk}

\author[J. Li]{Jiangpeng Li}
\address{School of Mathematics, Shandong University, Jinan 250100, China}
\email{jpli@mail.sdu.edu.cn}

\author[Y.-C. Sun]{Yu-Chen Sun}
\address{School of Mathematics, University of Bristol, Fry Building,
Woodland Road, Bristol BS8 1UG, UK}
\email{yuchensun93@163.com}

\author[H. Wang]{Hui Wang}
\address{The University of Hong Kong Shenzhen Institute of Research and Innovation,
Shenzhen, Guangdong 518057, China}
\email{wh0315@mail.sdu.edu.cn}

\makeatletter
\@ifundefined{subjclassname@2020}{%
  \@namedef{subjclassname@2020}{\textup{2020} Mathematics Subject Classification}}{}
\makeatother
\date{}
\hypersetup{pdftitle={On the lower bound for the L1 norm of Liouville exponential sums. II.},
pdfauthor={Ben Krause, Jiangpeng Li, Yu-Chen Sun, and Hui Wang}}
\begin{document}

\begin{abstract}
Let $\lambda$ denote the Liouville function, and put
\[
I_\lambda(X)=\int_0^1\left|\sum_{n\le X}\lambda(n)e^{2\pi i n\alpha}\right|\,d\alpha.
\]
In this paper, we prove that
\[
\limsup_{X \to \infty} \frac{\log I_{\lambda}(X)}{\log X} \ge \frac{1}{3}.
\]
\end{abstract}

\maketitle

\section{Introduction}

Let $\lambda$ be the completely multiplicative function defined by
$\lambda(p)=-1$ at every prime $p$. Write $e(t)=\exp(2\pi it)$, and equip
$\T=\R/\Z$ with Haar measure of total mass one. For $X\ge2$, put
\begin{equation}\label{eq:defI}
 S_X(\alpha)=\sum_{n\le X}\lambda(n)e(n\alpha),
 \qquad I_\lambda(X)=\int_\T|S_X(\alpha)|\,d\alpha.
\end{equation}
The first Fourier coefficient, Parseval's identity and the Cauchy-Schwarz
inequality give $1\le I_\lambda(X)\le X^{1/2}$.
The problem is to obtain lower
bounds that reflect the multiplicative structure of the coefficients.
Pandey and Radziwi\l\l~\cite[Theorem~1.1]{PR1} proved that
\begin{equation}\label{eq:quarter}
 I_\lambda(X)\gg_\eps X^{1/4-\eps}
\end{equation}
for every $\eps>0$. Their proof relates $I_\lambda(X)$ to zeros of
Dirichlet $L$-functions. More precisely, if $\chi$ is a fixed primitive
character and $\rho=\beta+i\gamma$ is a fixed zero of $L(s,\chi)$, then
\begin{equation}\label{eq:fixedzero}
 I_\lambda(X)\gg_{\rho,\chi,\eps}X^{\beta/2-\eps}
 \qquad(\eps>0)
\end{equation}
by~\cite[Theorem~1.2]{PR1}.
The $X^\eps$ loss in \eqref{eq:quarter} was removed
in~\cite[Theorem~1.1]{LSW}.

The true order of $I_\lambda(X)$ remains unknown. It is expected that
$I_\lambda(X)\gg X^{1/2-o(1)}$; see~\cite[\S1]{PR2}.
In this paper, we show that any upper bound $I_\lambda(X)\ll X^c$
valid for all sufficiently large $X$ must have $c\ge1/3$.

\begin{theorem}\label{thm:main}
Let
\[
 \ell=\liminf_{X\to\infty}\frac{\log I_\lambda(X)}{\log X},
 \qquad
 \calL=\limsup_{X\to\infty}\frac{\log I_\lambda(X)}{\log X}.
\]
Then $0\le\ell\le\calL\le1/2$ and
\[
 2\ell+\calL\ge1.
\]
In particular, $\calL\ge1/3$, and for every $\eps>0$ there are
arbitrarily large $N$ such that
\[
 I_\lambda(N)\ge N^{1/3-\eps}.
\]
\end{theorem}

\begin{remark}
The estimate \eqref{eq:quarter} holds at every sufficiently large scale.
Theorem~\ref{thm:main} gives the larger exponent $1/3$ in the limit
superior; it does not give $I_\lambda(X)\gg_\eps X^{1/3-\eps}$ for all
sufficiently large $X$.
\end{remark}

\noindent\textbf{Hypothesis \boldmath$\calH(\nu)$.}
Let $1/2\le\nu<1$. Every primitive Dirichlet $L$-function is zero-free in
$\Re s>\nu$.

Let us briefly discuss the motivation for the method and the key ideas.
Write $\T=\Maj\cup\Min$, where $\Maj$ and $\Min$ denote the major
and minor arcs, respectively.
The lower bounds in~\cite{PR1,LSW} come from information on the major arcs,
without using estimates on the minor arcs. Specifically, for a suitable
smooth weight $w$, we have
\[
I_{\lambda}(X) \ge \int_{\Maj}|\sum_{n}\lambda(n)w(n/X)e(n \alpha)| d \alpha.
\]
Here we use information
on the minor arcs.

We construct a kernel $G(\alpha)$ that resonates with $S_X(\alpha)$,
in the sense that
\[
 \int_\T S_X(\alpha)\overline{G(\alpha)}\,d\alpha\asymp X.
\]
Under $\calH(2/3)$, Proposition~\ref{prop:major} bounds the mean square
of $G(\alpha)$ over the major arcs using the hybrid large sieve for Dirichlet
polynomials. (If we assume the quasi-GRH bound in~\cite{OpenAI}, the
same argument applies with $\calH(7/8)$ in place of $\calH(2/3)$.)
With suitable choices of the scales, the major arcs
contribute $o(X)$ to this correlation, so the main term comes from
the minor arcs.

The kernel $G(\alpha)$ is chosen to have a bilinear structure.
More precisely, for a suitable parameter $Q\asymp X^u$ and
coefficients satisfying $|a_m| \ll m^{o(1)}$, we write
\[
G(\alpha)=\sum_{m \leq X/Q}a_m\sum_{n \leq Q} \lambda(n)e(mn \alpha).
\]
This construction is motivated by Vaughan's treatment of Waring's
problem~\cite[\S10]{Vaughan}, where restricting variables to smooth
numbers permits a bilinear decomposition of the exponential sums.
Here the product variables in $G(\alpha)$ provide this structure
directly. The large sieve inequality of $L^1$ type in
Lemma~\ref{lem:sampling} bounds $\sup_{\alpha \in \Min}|G(\alpha)|$ by $\frac{X^{1+\varepsilon}}{Q}I_{\lambda}(Q)$. Combining
this estimate with the resonance gives $X \ll \frac{X^{1+\varepsilon}}{Q}I_{\lambda}(Q)I_{\lambda}(X)$.

The rest of the paper is organized as follows. In Section 2, we construct the resonant kernel and deal with the minor arc case by using the bilinear structure. In Section 3, we collect mean-value theorems and a uniform bound for Dirichlet polynomials. In Section 4, we will handle the major arc cases and then prove our main theorem.

\subsection*{AI assistance}
The authors used ChatGPT (OpenAI) to assist with mathematical computations, the search for references and the preparation of the manuscript, but all the ideas come from the authors themselves. The authors take full responsibility for the results and their presentation.
\section{Minor arcs and bilinear structures}\label{sec:amplifier}

Fix $k\ge1$ and put $\kappa=2^{-1/(k+1)}$. For sufficiently large $Q,L$ and $Q \le X^{1/2-\varepsilon}$ with any $\varepsilon>0$, set
\begin{equation}\label{eq:scales}
 X=QL^k.
\end{equation}
Define
\begin{align}
 A_Q(\alpha)&=\sum_{\kappa Q<m\le Q}\lambda(m)e(m\alpha),
 \label{eq:AQ}\\
 G(\alpha)&=\sum_{\substack{\kappa Q<m\le Q\\
                     \kappa L<n_i\le L\ (1\le i\le k)}}
 \lambda(mn_1\cdots n_k)e(mn_1\cdots n_k\alpha).
 \label{eq:amplifier}
\end{align}
Every product in \eqref{eq:amplifier} lies in $(X/2,X]$. A simple computation gives
\begin{equation}\label{eq:correlation}
 \int_\T S_X(\alpha)\overline{G(\alpha)}\,d\alpha
 \sim(1-\kappa)^{k+1}X.
\end{equation}
Also,
\begin{equation}\label{eq:intervalnorm}
 \|A_Q\|_1\ll\log(2Q)I_\lambda(Q).
\end{equation}
Indeed, $A_Q=S_Q*K_Q$, where
$K_Q(\alpha)=\sum_{\kappa Q<n\le Q}e(n\alpha)$. Then \eqref{eq:intervalnorm}
follows from Young's inequality.

\begin{equation}\label{eq:majorarcs}
 \Maj_X(Q)=\bigcup_{q\le Q}\ \bigcup_{a\bmod q}^{*}
 \left\{\alpha\in\T:
 \left\|\alpha-\frac aq\right\|_\T\le\frac{Q}{qX}\right\},
 \qquad \Min_X(Q)=\T\setminus\Maj_X(Q).
\end{equation}

\begin{lemma}[{\cite[Lemma~2.1]{LSW}}]\label{lem:sampling}
Let $N\geqslant 1$ be an integer and $0<\delta\leqslant 1$.
Let $\mathcal{R}\subset\mathbb{T}$ be a finite set whose distinct
elements have circular distance at least $\delta$.
If $T\left(\alpha\right)=\sum_{|n|\leqslant N}b_n\mathrm{e}\left(n\alpha\right)$, then
\begin{equation}\label{eq:sampling}
 \sum_{r\in\mathcal{R}}|T\left(r\right)|
 \ll\left(N+\delta^{-1}\right)
 \int_{\mathbb{T}}|T\left(\alpha\right)|\mathrm{d}\alpha.
\end{equation}
\end{lemma}

\begin{lemma}\label{lem:minor}
If $Q^2\le X$, then for every $\eps>0$,
\begin{equation}\label{eq:minorbound}
 \sup_{\alpha\in\Min_X(Q)}|G(\alpha)|
 \ll_{k,\eps}\frac{X^{1+\eps}}{Q}I_\lambda(Q).
\end{equation}
\end{lemma}

\begin{proof}
Write $G(\alpha)$ as the following bilinear form:
\begin{equation}\label{eq:grouped}
 G(\alpha)=\sum_{j\le X/Q}b_k(j)A_Q(j\alpha),
\end{equation}
where $| b_k(j)|\ll_k j^{o(1)}$.

By
Dirichlet approximation, for any fixed $\alpha\in\Min_X(Q)$, there exists $a/q$ with $(a,q)=1$ such that
\[
 Q\le q\le \frac{X}{Q},\qquad
 \left|\alpha-\frac aq\right|\le\frac{Q}{qX}\le\frac1{q^2}.
\]
Split $j\le X/Q$ into blocks of at most $q/2$ consecutive integers.
If $h$ is the nonzero difference of two
indices in one block, then $|h|<q/2$ and
\[
 \|h\alpha\|_\T\ge\left\|\frac{ha}{q}\right\|_\T
 -\frac{|h|}{q^2}\ge\frac1{2q}.
\]
Since $Q\le q$, Lemma~\ref{lem:sampling}, applied to each block, gives
\begin{align*}
 \sum_{j\le X/Q}|A_Q(j\alpha)|
 &\ll\left(\frac{X}{Qq}+1\right)q\|A_Q\|_1\\
 &\ll\left(\frac XQ+q\right)\|A_Q\|_1
 \ll\frac XQ\|A_Q\|_1.
\end{align*}
Equations \eqref{eq:grouped} and \eqref{eq:intervalnorm} now give
\[
 |G(\alpha)|\ll_k X^{o(1)}\frac XQ\log(2Q)I_\lambda(Q)
 \ll_k\frac{X^{1+o(1)}}Q I_\lambda(Q).
\]
This proves \eqref{eq:minorbound}.
\end{proof}

\section{Mean-value estimates and a zero-free bound}\label{sec:analytic}

We collect the following two mean-value estimates, which will be used in the next section.

\begin{lemma}\label{lem:transfer}
Let $N\ge1$, and let $(b_n)$ be a sequence of complex numbers
supported in $(N/4,4N)$. Then, for all $X\ge1$,
\begin{equation}\label{eq:transfer}
 \int_{|\theta|\le1/X}\left|\sum_n b_ne(n\theta)\right|^2\,d\theta
 \ll\frac1N\int_\R\left|\sum_n b_nn^{-it}\right|^2
 \min\left(1,\frac{(N/X)^2}{1+|t|^2}\right)\,dt,
\end{equation}
with an absolute implied constant.
\end{lemma}

\begin{proof}
This is \cite[Lemma~4.1]{PR2}, with $b_n=a_nW(n/N)$ in their notation.
\end{proof}

\begin{lemma}\label{lem:hybrid}
Let $Q,T,N\ge1$. For complex coefficients $(a_n)_{n\le N}$
independent of the conductor and character,
\begin{equation}\label{eq:sieve}
 \sum_{q\le Q}\ \sideset{}{^*}\sum_{\chi\bmod q}
 \int_{-T}^T\left|\sum_{n\le N}a_n\chi(n)n^{-it}\right|^2\,dt
 \ll(Q^2T+N)\sum_{n\le N}|a_n|^2,
\end{equation}
where the starred sum runs over primitive characters modulo $q$.
\end{lemma}

\begin{proof}
This is the primitive hybrid large sieve \cite[Lemma~5.3]{PR2}.
\end{proof}

For $Y>0$, write $n\sim Y$ for $\kappa Y<n\le Y$.
For a primitive character $\chi\bmod q$ and an integer $D\ge1$,
define the Dirichlet polynomial
\begin{equation}\label{eq:FYD}
 F_{Y,D}(s,\chi)=
 \sum_{\substack{n\sim Y\\(n,D)=1}}
 \frac{\lambda(n)\chi(n)}{n^s}
 \qquad(s\in\mathbb C).
\end{equation}
The sum is empty when $Y<1$.

\begin{lemma}\label{lem:zerofree}
Assume $\calH(\vartheta)$ and fix $\vartheta<\sigma<1$. For every
$\eps>0$,
\begin{equation}\label{eq:zerofree}
 |F_{Y,D}(it,\chi)|
 \ll_{\vartheta,\sigma,\eps,k}
 Y^\sigma\{qD(3+|t|)\}^{\eps},
\end{equation}
uniformly in $q,D$, $Y\ge1$, and $t\in\R$.
\end{lemma}

\begin{proof}
For $\Re s>1$, define
\begin{equation}\label{eq:DS}
 \mathcal D_{\chi,D}(s)
 =\sum_{(n,D)=1}\frac{\lambda(n)\chi(n)}{n^s}
 =\frac{L(2s,\chi^2)}{L(s,\chi)}
  \prod_{p\mid D}(1+\chi(p)p^{-s}).
\end{equation}
The truncated Perron formula \cite[Corollary~5.3]{MV}
(cf.\ \cite[proof of Lemma~3.1]{PSW}) gives
\begin{equation}\label{eq:perron}
 F_{Y,D}(it,\chi)
 =\frac1{2\pi i}\int_{2-iY^2}^{2+iY^2}
 \mathcal D_{\chi,D}(s+it)
 \frac{Y^s-(\kappa Y)^s}{s}\,ds+O(1).
\end{equation}

Put $\sigma_0=(\vartheta+\sigma)/2$. Since the right-hand side of
\eqref{eq:DS} is holomorphic in $\Re s>\vartheta$, we shift the
contour in \eqref{eq:perron} to the line $\Re s=\sigma_0$.

For $u\ge\sigma_0$, \cite[Lemma~3]{Baker} and the argument in
\cite[\S14.2, pp.~336--337]{Titchmarsh} give, respectively,
\begin{equation}\label{eq:reciprocal}
 \begin{aligned}
 |L(u+iv,\chi)^{-1}|
 &\ll_{\vartheta,\sigma}
 \{q(3+|v|)\}^{o(1)}\qquad(q>1),\\
 |\zeta(u+iv)^{-1}|
 &\ll_{\vartheta,\sigma}(3+|v|)^{o(1)}.
 \end{aligned}
\end{equation}
Since $\sigma_0>1/2$,
$|L(2u+2iv,\chi^2)|\le\zeta(2\sigma_0)$ and
$\prod_{p\mid D}(1+p^{-\sigma_0})\ll_{\sigma_0}D^{o(1)}$.
Hence
\begin{equation}\label{eq:DSbound}
 |\mathcal D_{\chi,D}(u+iv)|
 \ll_{\vartheta,\sigma}
 \{qD(3+|v|)\}^{o(1)}\qquad(u\ge\sigma_0).
\end{equation}
\eqref{eq:zerofree} follows by combining everything.
\end{proof}

\section{The major arc estimate and the proof of Theorem \ref{thm:main}}\label{sec:major}

\begin{proposition}\label{prop:major}
Assume $\calH(\vartheta)$, and fix $\vartheta<\sigma<1$ and $k\ge1$.
Let $Q,L \to \infty$ as $X \to \infty$, and suppose that \eqref{eq:scales} holds with $Q^2\le X$. For every
$\eps>0$,
\begin{equation}\label{eq:majorbound}
 \int_{\Maj_X(Q)}|G(\alpha)|^2\,d\alpha
 \ll_{\vartheta,\sigma,k,\eps}
 X^\eps\left(Q^2L^{2\sigma-1}
       +X\left(\frac{Q^2}{X}\right)^{2(1-\sigma)/k}\right).
\end{equation}
\end{proposition}

\begin{proof}
Throughout the proof, $\vartheta,\sigma,k$ are fixed. Recall that
\[
 G(\alpha)=\sum_{\substack{m\sim Q\\n_1,\ldots,n_k\sim L}}
 \lambda(mn_1\cdots n_k)e(mn_1\cdots n_k\alpha).
\]
Let
\[
 \mathcal G_{a,q}(t)
 =\sum_{\substack{m\sim Q\\n_1,\ldots,n_k\sim L}}
 \frac{\lambda(mn_1\cdots n_k)}{(mn_1\cdots n_k)^{it}}
 e\!\left(\frac{amn_1\cdots n_k}{q}\right).
\]
By the definition of $\Maj_X(Q)$ and Lemma~\ref{lem:transfer}, we have
\begin{align}
\int_{\Maj_X(Q)}|G(\alpha)|^2\,d\alpha
 & \le\sum_{q\le Q}\ \sideset{}{^*}\sum_{a\bmod q}
 \int_{-Q/(qX)}^{Q/(qX)}
 \left|G\left(\frac aq+\beta\right)\right|^2\,d\beta \notag \\
 & \ll \frac1X\sum_{q\le Q}\ \sideset{}{^*}\sum_{a\bmod q}\int_\R|\mathcal G_{a,q}(t)|^2
 \min\left(1,\frac{(Q/q)^2}{1+|t|^2}\right)\,dt. \label{eq:energy}
\end{align}
Split the sum in $\mathcal G_{a,q}$ according to
\[
 e_0=(m,q),\qquad
 e_i=\left(n_i,\frac{q}{e_0\cdots e_{i-1}}\right)
 \quad(1\le i\le k).
\]
Writing $m=e_0u_0$, $n_i=e_iu_i$, $\vec e=(e_0,\ldots,e_k)$
and $E=e_0\cdots e_k$, we obtain
\begin{equation}\label{eq:reducedscales}
 \mathcal G_{a,q}(t)
 =\sum_{\substack{\vec e\\E\mid q}}
 \frac{\lambda(E)}{E^{it}}
 \sum_{\substack{u_0\sim Q/e_0,\ u_i\sim L/e_i\ (1\le i\le k)\\
 (u_i,q/(e_0\cdots e_i))=1\ (0\le i\le k)}}
 \frac{\lambda(u_0\cdots u_k)}{(u_0\cdots u_k)^{it}}
 e\!\left(\frac{aEu_0\cdots u_k}{q}\right).
\end{equation}
For each $\vec e$, put
\[
 r=q/E,\qquad Y_0=Q/e_0,\qquad Y_i=L/e_i\quad(1\le i\le k).
\]
Since $(u_i,r)=1$, we insert the character identity
\[
 e(b/r)=\frac1{\varphi(r)}
 \sum_{\chi\bmod r}\tau_r(\overline\chi)\chi(b)
 \qquad((b,r)=1),
\]
where $\tau_r(\chi)=\sum_{b\bmod r}\chi(b)e(b/r)$, into
\eqref{eq:reducedscales}. There are at most
$d_{k+2}(q)\ll q^{o(1)}$ vectors $\vec e$.
The Cauchy-Schwarz inequality over
$\vec e$ gives
\[
|\mathcal G_{a,q}(t)|^2
 \leq q^{o(1)}
 \sum_{\substack{\vec e\\E\mid q}}\frac1{\varphi(r)^2}\times\left|
 \sum_{\chi\bmod r}\tau_r(\overline\chi)\chi(a)
 \prod_{i=0}^k
 \sum_{\substack{u\sim Y_i\\(u,q/(e_0\cdots e_i))=1}}
 \frac{\lambda(u)\chi(u)}{u^{it}}\right|^2.
\]
 
Fix $\vec e$. Since $r\mid q$, any $a \in \Z_q^{*}$ can be written as $sr+b \frac{q}{(q,r)}$ with $s \in \Z_{\frac{q}{(q,r)}}^{*}$ and $b \in \Z_{r}^{*}$. Thus, for $\chi_1,\chi_2\bmod r$,
\[
 \begin{aligned}
 \sum_{\substack{a\bmod q\\(a,q)=1}}
\chi_1(a)\overline{\chi_2(a)} &=\sum_{\substack{s\bmod \frac{q}{(q,r)}\\(s,\frac{q}{(q,r)})=1}}
 \sum_{\substack{b\bmod r\\(b,r)=1}}\chi_1(b)\overline{\chi_2(b)}\\
 &=\frac{\varphi(q)}{\varphi(r)}
 \sum_{\substack{b\bmod r\\(b,r)=1}}
 \chi_1(b)\overline{\chi_2(b)}\\
 &=\varphi(q)\mathbf1_{\chi_1=\chi_2}.
 \end{aligned}
\]
Sum the Cauchy-Schwarz bound over reduced $a\bmod q$, expand the square,
and apply this identity. We obtain
\begin{equation}\label{eq:characterenergy}
\sideset{}{^*}\sum_{a\bmod q}|\mathcal G_{a,q}(t)|^2
 \ll_k q^{o(1)}
 \sum_{\substack{\vec e\\E\mid q}}
 \frac{\varphi(q)}{\varphi(r)^2}
 \sum_{\chi\bmod r}|\tau_r(\overline\chi)|^2
 \times
 \left|\prod_{i=0}^k
 \sum_{\substack{u\sim Y_i\\(u,q/(e_0\cdots e_i))=1}}
 \frac{\lambda(u)\chi(u)}{u^{it}}\right|^2.    
\end{equation}

Let $\chi^*\bmod r_0$ be the primitive character inducing
$\chi\bmod r$, and write $r=r_0r_1$. This means that
$\chi(n)=\chi^*(n)$ for $(n,r)=1$, and $\chi(n)=0$ otherwise.
Since $\chi^*(n)=0$ when $(n,r_0)>1$, M\"obius inversion gives
\[
 \chi(n)=\chi^*(n)\mathbf1_{(n,r_1)=1}
 =\chi^*(n)\sum_{d\mid(n,r_1)}\mu(d).
\]
Consequently,
\begin{equation}\label{eq:induced}
 \begin{aligned}
 \tau_{r_0r_1}(\chi)
 &=\sum_{d\mid r_1}\mu(d)\chi^*(d)
   \sum_{b\bmod r_0r_1/d}\chi^*(b)e\!\left(\frac{db}{r_0r_1}\right)\\
 &=\mu(r_1)\chi^*(r_1)\tau_{r_0}(\chi^*).
 \end{aligned}
\end{equation}
By \cite[Theorem~9.7]{MV},
$|\tau_{r_0}(\overline{\chi^*})|^2=r_0$. Together with
\eqref{eq:induced}, this gives
\begin{equation}\label{eq:Eoverh}
 \frac{\varphi(q)|\tau_r(\overline\chi)|^2}{\varphi(r)^2}
 \le\frac{\varphi(Er_0r_1)r_0}{\varphi(r_0r_1)^2}
 \le\frac{E}{r_1}\left(\frac{r_0r_1}{\varphi(r_0r_1)}\right)^2
 \ll X^{o(1)}\frac{E}{r_1}.
\end{equation}
Put $D_i=r_1e_{i+1}\cdots e_k$.
The short sums in \eqref{eq:characterenergy} satisfy
\begin{equation}\label{eq:conductorcoprimality}
 \sum_{\substack{u\sim Y_i\\(u,q/(e_0\cdots e_i))=1}}
 \frac{\lambda(u)\chi(u)}{u^{it}}
 =\sum_{\substack{u\sim Y_i\\(u,r_0D_i)=1}}
 \frac{\lambda(u)\chi^*(u)}{u^{it}}
 =F_{Y_i,D_i}(it,\chi^*).
\end{equation}
Substituting these identities and \eqref{eq:Eoverh} into
\eqref{eq:characterenergy}, then into \eqref{eq:energy}, and
writing $q=Er_0r_1$, we obtain
\begin{align}
 \int_{\Maj_X(Q)}|G(\alpha)|^2\,d\alpha
 &\ll_k X^{o(1)}
 \sum_{\substack{\vec e,r_1\\Er_1\le Q}}\frac E{r_1X}
 \sum_{r_0\le Q/(Er_1)}\ \sideset{}{^*}\sum_{\chi\bmod r_0}\notag\\
 &\qquad\times
 \int_\R\left|\prod_{i=0}^kF_{Y_i,D_i}(it,\chi)\right|^2
 \min\left(1,\frac{(Q/(Er_0r_1))^2}{1+|t|^2}\right)\,dt.
 \label{eq:primitiveenergy}
\end{align}

Fix $\vec e,r_1$ and choose $j\in\{1,\ldots,k\}$ with
$Y_j$ maximal. Since $E\le Q$,
\begin{equation}\label{eq:survivor}
 \prod_{i=1}^kY_i=\frac{L^k}{e_1\cdots e_k}\ge\frac{L^k}{Q},
 \qquad \frac{L}{Q^{1/k}}\le Y_j\le L.
\end{equation}
Write
\[
 \prod_{\substack{0\le i\le k\\i\ne j}}F_{Y_i,D_i}(it,\chi)
 =\sum_{n\le N}a_n\chi(n)n^{-it},
 \qquad N=\frac{X}{EY_j}.
\]
The coefficients are independent of $r_0,\chi$ and satisfy
$|a_n|\le d_k(n)$, so $\sum_{n\le N}|a_n|^2\ll_k NX^{o(1)}$.

Split $r_0\le Q/(Er_1)$ into dyadic intervals
$R_0\le r_0<2R_0$, and put $T=Q/(Er_1R_0)$.
For $r_0\le2R_0$ and $|t|\le2^{\ell+1}T$, we have
$r_0D_j(3+|t|)\ll X2^\ell$. Lemma~\ref{lem:zerofree} therefore gives
\[
 |F_{Y_j,D_j}(it,\chi)|^2
 \ll_{\vartheta,\sigma,k}Y_j^{2\sigma}X^{o(1)}2^{\ell/2}
 \qquad(\ell\ge0).
\]

Splitting the $t$-integral into $|t|\le2T$ and
$2^\ell T<|t|\le2^{\ell+1}T$ for $\ell\ge1$ and applying Lemma~\ref{lem:hybrid}, we obtain
\begin{align*}
 &\sum_{R_0\le r_0<2R_0}\ \sideset{}{^*}\sum_{\chi\bmod r_0}\int_\R
 \left|\prod_{i=0}^kF_{Y_i,D_i}(it,\chi)\right|^2
 \min\left(1,\frac{(Q/(Er_0r_1))^2}{1+|t|^2}\right)\,dt\\
 &\quad\ll X^{o(1)}Y_j^{2\sigma}
 \sum_{\ell\ge0}2^{-3\ell/2}
 \sum_{r_0\le2R_0}\ \sideset{}{^*}\sum_{\chi\bmod r_0}
 \int_{-2^{\ell+1}T}^{2^{\ell+1}T}
 \left|\sum_{n\le N}a_n\chi(n)n^{-it}\right|^2\,dt\\
 &\quad\ll X^{o(1)}Y_j^{2\sigma}
 \sum_{\ell\ge0}2^{-3\ell/2}(2^\ell R_0^2T+N)N\\
 &\quad\ll X^{o(1)}Y_j^{2\sigma}
 \left(\frac{Q^2}{E^2r_1^2}+N\right)N.
\end{align*}
The height sums converge, and
$R_0^2T=QR_0/(Er_1)\le Q^2/(E^2r_1^2)$.
Summing over the $O(\log(2Q))$ conductor intervals and substituting
$N=X/(EY_j)$, we obtain
\begin{align}
 &\frac E{r_1X}
 \sum_{r_0\le Q/(Er_1)}\ \sideset{}{^*}\sum_{\chi\bmod r_0}\int_\R
 \left|\prod_{i=0}^kF_{Y_i,D_i}(it,\chi)\right|^2
 \min\left(1,\frac{(Q/(Er_0r_1))^2}{1+|t|^2}\right)\,dt\notag\\
 &\qquad\ll X^{o(1)}
 \left(\frac{Q^2Y_j^{2\sigma-1}}{E^2r_1^3}
       +\frac{XY_j^{2\sigma-2}}{Er_1}\right).
 \label{eq:twoterms}
\end{align}
Since $1/2<\sigma<1$, \eqref{eq:survivor} implies
\[
 Y_j^{2\sigma-1}\le L^{2\sigma-1},\qquad
 Y_j^{2\sigma-2}\le\left(\frac{Q^2}{X}\right)^{2(1-\sigma)/k}.
\]
Combining the preceding estimates, we obtain \eqref{eq:majorbound}.
\end{proof}

\begin{proof}[Proof of Theorem~\ref{thm:main}]
It suffices to prove $2\ell+\calL\ge1$.

If some primitive Dirichlet $L$-function has a zero
$\rho=\beta+i\gamma$ with $\beta>2/3$, then
\eqref{eq:fixedzero} implies
\[
 \ell\ge\frac\beta2>\frac13,
\]
so $2\ell+\calL\ge3\ell>1$. We may therefore assume
$\calH(2/3)$.

Fix $k\ge2$. For $X\to\infty$, put
\begin{equation}\label{eq:powerchoice}
 Q=X^{(k-1)/(2k-1)},\qquad L=X^{1/(2k-1)}.
\end{equation}
Then $X=QL^k$, $Q^2<X$, and $Q,L\to\infty$.
Proposition~\ref{prop:major}, applied with $\vartheta=2/3$ and
$\sigma=3/4$, gives
\begin{equation}\label{eq:smallenergy}
 \int_{\Maj_X(Q)}|G(\alpha)|^2\,d\alpha
 \ll_k X^{1+o(1)}
 \left(L^{-1/2}+L^{-1/(2k)}\right)=o(X).
\end{equation}
The Cauchy-Schwarz inequality and $\|S_X\|_2^2=X$ give
\[
 \left|\int_{\Maj_X(Q)}S_X(\alpha)\overline{G(\alpha)}\,d\alpha\right|
 \le\|S_X\|_2
 \left(\int_{\Maj_X(Q)}|G(\alpha)|^2\,d\alpha\right)^{1/2}
 =o(X).
\]
Subtracting this from \eqref{eq:correlation}, we obtain
\[
 X\ll_k\left|\int_{\Min_X(Q)}S_X(\alpha)\overline{G(\alpha)}\,d\alpha\right|
 \le I_\lambda(X)\sup_{\alpha\in\Min_X(Q)}|G(\alpha)|.
\]
Lemma~\ref{lem:minor} therefore yields
\begin{equation}\label{eq:twoscale}
 I_\lambda(X)I_\lambda(Q)
 \gg_k QX^{-o(1)}.
\end{equation}
Our result follows by letting $k\to\infty$.
\end{proof}

\end{document}